\documentclass[11pt]{article}

\usepackage{amsmath}
\usepackage{amssymb}
\usepackage{latexsym}

\newtheorem{assumption}{Assumption}

\def\qed{ \ \vrule width.2cm height.2cm depth0cm\smallskip}
\def \ep{\hbox{ }\hfill$\Box$}

\newenvironment{proof}{\noindent {\bf Proof.\/}}{$\qed$\vskip 0.1in}

\newcommand{\la}{\langle}
\newcommand{\ra}{\rangle}

\newcommand{\ba}{\begin{array}}
\newcommand{\ea}{\end{array}}
\newcommand{\be}{\begin{equation}}
\newcommand{\ee}{\end{equation}}
\newcommand{\bea}{\begin{eqnarray}}
\newcommand{\eea}{\end{eqnarray}}
\newcommand{\beaa}{\begin{eqnarray*}}
\newcommand{\eeaa}{\end{eqnarray*}}

\def\cD{{\cal D}}

\def\cF{{\cal F}}

\def\cH{{\cal H}}
\def\cI{{\cal I}}

\def\cL{{\cal L}}
\def\cM{{\cal M}}

\def\cP{{\cal P}}

\def\cS{{\cal S}}

\def\dbE{\mathbb{E}}
\def\dbF{\mathbb{F}}

\def\dbH{\mathbb{H}}
\def\dbI{\mathbb{I}}

\def\dbL{\mathbb{L}}

\def\dbP{\mathbb{P}}
\def\dbR{\mathbb{R}}
\def\dbS{\mathbb{S}}

\def\hL2{\widehat{\dbL}^2}
\def\hH2{\widehat{\cH}^2}
\def\hD2{\widehat{\cD}^2}
\def\hI2{\widehat{\cI}^2}
\def\hS2{\widehat{\cS}^2}

\def\a{\alpha}

\def\g{\gamma}
\def\d{\delta}
\def\e{\varepsilon}

\def\l{\lambda}

\def\t{\tau}
\def\f{\varphi}
\def\th{\theta}
\def\o{\omega}

\def\O{\Omega}
\def\q{\quad}

\def\pa{\partial}
\def\cd{\cdot}
\def\cds{\cdots}

\def\tr{\hbox{\rm tr}}

\def\qed{ \hfill \vrule width.25cm height.25cm depth0cm\smallskip}

\newcommand{\basa}{\begin{assumption}}
\newcommand{\easa}{\end{assumption}}

\newcommand{\bas}{\begin{assum}}
\newcommand{\eas}{\end{assum}}

\def\limsup{\mathop{\overline{\rm lim}}}

\def\esup{\mathop{\rm ess\;sup}}
\def\einf{\mathop{\rm ess\;inf}}

\def\pa{\partial}

 \def\cd{\cdot}
\def\cds{\cdots}

\def\tr{\hbox{\rm tr$\,$}}

\def\dis{\displaystyle}

\def\cad{{c\`{a}dl\`{a}g}}

\def\1{\mathbf{1}}

\def\:{\!:\!}
\def\reff#1{{\rm(\ref{#1})}}
\def \proof{{\noindent \it Proof.\quad}}

at 9pt

\begin{document}

\newtheorem{thm}{Theorem}[section]
\newtheorem{lem}[thm]{Lemma}
\newtheorem{cor}[thm]{Corollary}
\newtheorem{prop}[thm]{Proposition}
\newtheorem{rem}[thm]{Remark}
\newtheorem{eg}[thm]{Example}
\newtheorem{defn}[thm]{Definition}
\newtheorem{assum}[thm]{Assumption}
\newtheorem{example}[thm]{Example}

\renewcommand {\theequation}{\arabic{section}.\arabic{equation}}
\def\thesection{\arabic{section}}

\title{Martingale Representation Theorem for the $G$-expectation
\footnote{Authors would like to thank the anonymous referees and Marcel Nutz for careful
reading of the first draft and numerous useful comments. }}
\author{H. Mete {\sc Soner}\footnote{ETH (Swiss Federal Institute of Technology),
Zurich, hmsoner@ethz.ch and Swiss Finance Institute. Research partly supported by the
European Research Council under the grant 228053-FiRM.
Financial support from
the ETH Foundation
is also gratefully acknowledged.}
       \and Nizar {\sc Touzi}\footnote{CMAP, Ecole Polytechnique Paris, nizar.touzi@polytechnique.edu.
       Research supported by the Chair {\it Financial Risks} of the {\it Risk Foundation} sponsored by Soci\'et\'e
              G\'en\'erale, the Chair {\it Derivatives of the Future} sponsored by the {F\'ed\'eration Bancaire Fran\c{c}aise}, and
              the Chair {\it Finance and Sustainable Development} sponsored by EDF and Calyon. }
       \and Jianfeng {\sc Zhang}\footnote{University of Southern California, Department of Mathematics,
       jianfenz@usc.edu. Research supported in part by NSF grant DMS 06-31366.}
}
\date{First version: January 20, 2010\\This version: September 5,  2010}

\maketitle

\begin{abstract}
This paper considers the  nonlinear theory of $G$-martingales
as introduced by Peng in \cite{Peng1,Peng}.
A martingale representation theorem for
this theory is proved by using the techniques and the results
established in \cite{STZ09} for the second order stochastic target
problems and the second order backward stochastic differential equations.
In particular, this representation provides a hedging strategy
in a market with an uncertain volatility.

\vspace{5mm}


\end{abstract}
\newpage

\section{Introduction}
\setcounter{equation}{0}

The notion of a  $G$-expectation as recently introduced by Peng  \cite{Peng1,Peng}
has several motivations and applications.
One of them is the study of financial problems with uncertainty about the volatility.
This  important problem was also considered earlier by Denis \& Martini \cite{DM}.
Motivated by this application,
Denis \& Martini developed an almost pathwise theory of stochastic calculus.
In this second approach,  probabilistic statements are required to hold  {\it{quasi surely}}: namely
$\dbP$-almost surely for all probability measures $\dbP$ from a large
class of mutually singular measures $\cP$.
Denis \& Martini employ functional analytic techniques while Peng's approach utilizes the
theory of viscosity solutions of parabolic partial differential equations.

Indeed, the $G$-expectation is defined by Peng using the nonlinear
heat equation,
\beaa
-\pa_t u - G(D^2 u) = 0 ~~\mbox{on}~[0,1),
\eeaa
where the time maturity is taken to be $T=1$ and  for given
$d\times d$ symmetric matrices $ \overline{a} >0 $ and
$0 \le \underline{a}\le\overline{a}$, the nonlinearity
$G$ is defined by, 
\be
\label{e.g}
G(\gamma):=\frac12\sup\{\ \tr[\gamma a] \ | \  \underline{a}\le a\le \overline{a}\},\q \g\in \dbR^{d\times d}.
\ee
Then for  ``Markov-like" random variables,
the $G$-expectation and conditional
expectations are defined through the solution of the above equation
with this random variable as its terminal condition at time $T=1$.
A $G$-martingale is
then defined easily as a process which satisfies the martingale property
by this  conditional expectation.  A brief introduction
to this theory is provided in Section \ref{s.G} below.

Denis \& Martini \cite{DM} also construct a
similar structure of quasi-sure stochastic analysis.
However, they use a quite
different approach which utilizes the set $\cP$ of all probability
measures $\dbP$ so that the canonical map in the Wiener space
is a martingale under $\dbP$ and the quadratic variation of this
martingale lies between $\underline{a}\le\overline{a}$.
Although the constructions of the quasi sure analysis and the
$G$-expectations are substantially different, these theories
are very closely related as proved recently by
Denis, Hu \& Peng \cite{DHP}.  The paper
\cite{DHP} also provides a dual representation
of the $G$-expectation as the supremum of
expectations  over $\cP$.  This duality and
more generally  the dynamic programming
principle is generalized by Nutz \cite{nutz}
who considers lower and upper bounds 
$\underline{a}$, $\overline{a}$ that are random processes. 

A probabilistic construction similar to quasi-sure stochastic analysis and  $G$-expectations,
is the  theory of second order backward stochastic
differential equations (2BSDE).  This theory is developed in \cite{CST,cstv,STgammaSIAM} as a generalization
of BSDEs as initially introduced in \cite{EPQ,pardouxpeng}.
In particular, 2BSDEs provide a stochastic representation for fully nonlinear partial differential
equations.  Since the $G$-expectation is defined through such a nonlinear equation, one
expects that the $G$-expectations are naturally connected to the
2BSDEs.  Equivalently, 2BSDEs can be viewed as the extension of $G$-expectations
to more general nonlinearities.  Indeed, recently the authors developed
such a generalization and a duality theory for 2BSDEs using probabilistic
constructions similar to quasi-sure analysis \cite{STZ09a,STZ09,STZ09b}.

In this paper, we investigate the question of representing
an arbitrary $G$-martingale in terms of stochastic integrals
and other processes.
Specifically, we fix a finite horizon say $T=1$.
Since all martingales can be seen as the conditional
expectation, we also fix  the final value $\xi$.
We then would like to construct stochastic processes
$H$ and $K$ so that
$$
Y_t := \dbE^G_t[\xi] = \xi -\int_t^1 \ H_s dB_s +K_1 - K_t\
= \dbE^G[\xi] + \int_0^t \ H_s dB_s - K_t,
$$
where $\dbE^G_t$ is the $G$-conditional expectation
and the process $M:=-K$ is a non-increasing $G$-martingale.
The stochastic integral that appears in the above is  the
regular It\^o one.  But it is also defined quasi-surely.
More precisely, the above statement holds almost-surely for all
probability measures in $\cP$.  Equivalently, the above equation holds
quasi-surely in the sense of Denis \& Martini.
In particular, all the above processes
as well as the stochastic integral are defined on the support of all
measures in the set $\cP$.  This is an important
property of this martingale representation
as  $\cP$ contains
measures which are mutually singular.
Moreover,  there is no
measure that dominates all measures in $\cP$.
Hence the above processes are defined
on a large subset of our probability space.

A partial answer to this question
was already provided by Xu and Zhang \cite{XZ} for the
class of {\em{symmetric $G$-martingales}},
i.e. a process $N$ which is
both itself and  $-N$
are $G$-martingales.  Since the $G$-expectation is
not  linear, the class of symmetric martingales
is a strict subset of all $G$-martingales.
In particular, the representation of
symmetric martingales are obtained using
only the stochastic integrals.
We obtain the martingale representation
in Theorem \ref{t.main} for almost all square-integrable
martingales.
This result essentially provides a complete answer to the question of representation
for the integrable classes defined in \cite{Peng}.

Our analysis utilizes the already mentioned
duality result of Denis, Hu and Peng \cite{DHP}.
Similar to \cite{DHP}, we also provide a dual characterization of $G$-martingales
as an immediate consequence of the results in \cite{DHP,Peng}.
This observation is one of the key-ingredients of our representation proof.
Moreover, it can be used to extend the definition of $G$-martingales to a
class larger than the integrability class $\cL^1_G$ of Peng.
Indeed, the above martingale representation result could also be proved
for a larger class of random variables.  But this
development also requires the extension of $G$-expectations
and conditional expectations to this larger class.  These types of results
are not pursued  here.  But in an example, Example \ref{e.l1} below, we show that
the integrability class $\cL^1_G$  does not include all
bounded random variables.  Thus it is desirable to extend the theory to
a larger class of random variables using the equivalent definitions that do not
refer to partial differential equations.  Indeed such a theory is developed
by the authors in \cite{STZ09a,STZ09,STZ09b}.

The paper is organized as follows. In Section \ref{s.G}, we review
the theory of $G$-expectations and $G$-martingales.
Section \ref{s.qs} defines the quasi-sure analysis of Denis \& Martini
and also provides the dual formulation. The main
ingredients for our approach, such as the norms and spaces,
are collected in Section \ref{s.norms}.
The main result is then stated and proved in
Section \ref{s.mrt}.
In the Appendix, we provide an approximation argument
for the solutions of the partial differential equation.
Then the connection between the integrability class
of Peng and the spaces utilized in this paper is given
in the subsection \ref{ss.lp}.

After the completion and the submission of this manuscript, 
we became aware of the manuscript of Song
\cite{song} which proves a decomposition result 
for random variables in $\cL^p_G$ with $p >1$.  
He obtained this result
after a preliminary version of this manuscript,
without Lemma \ref{l.equal}, below, was circulated. 
Indeed, it is clear that  a slight 
 extension of Theorem \ref{t.main},  below, to $\cL^p_\cP$,
 together with Lemma \ref{l.equal},
 implies the decomposition result 
 \reff{e.main} for any $\xi \in \cL^p_\cP$ with $p>1$.
  We also emphasize that, in contrast with \cite{song}, 
this manuscript
 considers the possibly degenerate case $\underline{a}\ge 0$,
 see Assumption \ref{e.matrix}.

\subsection{Notation and spaces}
\label{ss.notations}
We collect all the spaces and the notation used in the paper with a reference
to their definitions. We always assume that $ \overline{a} >0$, $ 0\le  \underline{a}\le \overline{a}$.\\
$\bullet$\hspace{5pt} $\dbF =\{\cF^B_{t}, t\ge 0\}$ is the filtration generated by the
canonical process  $B$.\\
$\bullet$\hspace{5pt} $\dbE^G$ is the $G$-expectation, defined in \cite{Peng}
and in subsection \ref{ss.g}.\\
$\bullet$\hspace{5pt} $\dbE^G_t$ is the conditional $G$-expectation.\\
$\bullet$\hspace{5pt} $\cL_{ip}$ is the space of random variables of the form $\f(B_{t_1},\cds, B_{t_n})$ with
a bounded, Lipschitz deterministic function $\f$ and
time points $0 \le t_1 \le \ldots \le t_n \le 1$.\\
$\bullet$\hspace{5pt} $\cL^p_G$ is the integrability class defined in subsection \ref{ss.g}
as the closure of $\cL_{ip}$.\\
$\bullet$\hspace{5pt} $\cH^{p,0}_G$ is the space of
piecewise constant $G$-stochastic integrands, see subsection \ref{ss.Gintegral}.\\
$\bullet$\hspace{5pt} $\cH^p_G$ is the integrability class defined in subsection \ref{ss.Gintegral} as the closure of $\cH^{p,0}_G$.\\
$\bullet$\hspace{5pt}  $\cP=\overline\cP^W_{[\underline a, \overline a]}$
measures under which the canonical process is a martingale and satisfies \reff{abound}.\\
$\bullet$\hspace{5pt} $\cP(t,\dbP)$ is defined in \reff{e.ptp}.\\
$\bullet$\hspace{5pt}  $\dbL^p_{\cP}$ is the set of all $p$-integrable random variables; see \reff{e.lnorm}. \\
$\bullet$\hspace{5pt} $\cL^p_\cP$ is the the closure of $\cL_{ip}$ under
the norm $\dbL^p_{\cP}$; see \reff{e.lnorm}.\\
$\bullet$\hspace{5pt} $\dbH^p_{\cP}$  is the set of all $p$-integrable, $\dbR^d$-valued
stochastic integrands; see \reff{e.hnorm}.\\
$\bullet$\hspace{5pt} $\cH^p_\cP$ is the closure of ${\cH^{p,0}_G}$
under the norm $\| \cdot \|_{\dbH^p_\cP}$; see Definition \ref{d.hnorm}\\
$\bullet$\hspace{5pt} $\dbS^p_\cP$  is the set of all $p$-integrable, {\it continuous}
processes; see Definition \ref{d.hnorm}.\\
$\bullet$\hspace{5pt} $\dbI^p_\cP$  is the subset of $\dbS^p_\cP$  that
are {\it non-decreasing} with initial value $0$; see Definition \ref{d.hnorm}.\\
$\bullet$\hspace{5pt} $\dbS_d$ is the set of all $d\times d$ symmetric matrices with the usual ordering and identity $I_d$.   \\
$\bullet$\hspace{5pt}  For $\nu, \eta \in \dbR^d$, $A:=\nu\otimes \eta\in \dbS_d $ is defined by $Ax= (\eta\cdot x) \nu$
for any $x \in \dbR^d$.\\
$\bullet$\hspace{5pt} For $A\in \dbS_d$, $\nu_k \in \dbR^d$ are its orthonormal
eigenvectors and $\lambda_k$ are the corresponding eigenvalues so that
$$
A = \sum_k\ \lambda_k [\nu_k \otimes \nu_k].
$$ \\
$\bullet$\hspace{5pt} For $A\in \dbS_d$, and a real number, $A \vee cI_d \in \dbS_d$ is defined by
$$
A \vee cI_d := \sum_k\ (\lambda_k \vee c ) \  [\nu_k \otimes \nu_k].
$$

\section{$G$-stochastic analysis of Peng \cite{Peng1,Peng}}
\label{s.G}
\setcounter{equation}{0}

We fix the time horizon $T=1$.  Let $\O:= \{\o\in C([0,1], \dbR^d): \o(0) = 0\}$ be 
the canonical space, $B$ the canonical process, and $\dbP_0$ the 
Wiener measure. $\dbF =\{\cF^B_t, t\in [0,1]\}$ is the filtration generated
by $B$. We note that $\cF^B_{t-} = \cF^B_t \neq \cF^B_{t+}$.

In what follows, we always use the space $\O$ together with the filtration 
$\dbF$. We remark that we do not augment the filtration, as usually done 
in standard stochastic analysis literature. In fact, for any probability measure 
$\dbP$ on $(\O, \cF_1)$, denote by $\bar\dbF^\dbP = \{\bar\cF^\dbP_t, 0\le t\le 1\}$ 
the augmented filtration of $\dbF$ under $\dbP$, we have the following straightforward result.

\begin{lem}
\label{lem-filtration}
For any $\bar\cF^\dbP_t$-measurable random variable $\xi$, 
there exists a unique {\rm{(}}$\dbP$-a.s.{\rm{)}}
$\cF_t$-measurable random variable $\tilde \xi$ such that $\tilde \xi = \xi$, 
$\dbP$-a.s.. 

Similarly, for every $\bar\dbF^\dbP$-progressively measurable process
$X$, there exists a unique  
$\dbF$-progressively measurable process $\tilde X$ such that $\tilde X = X$, 
$dt\times d \dbP$-a.s.. Moreover, if $X$ is $\dbP$-almost surely continuous, 
then one can choose $\tilde X$ to be 
$\dbP$-almost surely continuous.
\end{lem}

\proof  Lemma 2.4 in \cite{STZ09a} proves 
the analogous result for the right continuous 
filtration $\dbF^+:=\{\cF^B_{t+}, 0\le t\le 1\}$ and its augmentation, 
instead of  $\dbF$ and its augmentation.
However, the proof does not change in this context
and we prove the above result
following the proof Lemma 2.4 in \cite{STZ09a} line by line.
\ep

In what follows, quite often we make use of the above
result.  Indeed,
when a probability measure $\dbP$ is given, 
we will consider any process in its $\dbF$-progressively measurable version. 
However, we emphasize that these versions, in general,
may depend on $\dbP$.

\subsection{$G$-expectation and $G$-martingale}
\label{ss.g}
Following Peng \cite{Peng1}, 
let $G$ be as in \reff{e.g} with two given
$d\times d$ symmetric matrices satisfying
\begin{equation}
\label{e.matrix}
0\le  \underline{a}\le \overline{a},
\qquad
\overline{a} >0.
\end{equation}
Notice that we allow
degenerate 
diffusion matrices as the only
positivity assumption 
is placed on the upper bound.

For a bounded  Lipschitz continuous
function $\f$ on $\dbR^d$, let $u$ be the
unique,  bounded, Lipschitz continuous viscosity solution 
of the following parabolic equation,
\bea
\label{Gpde1}
-\pa_t u - G(D^2 u) = 0 ~~\mbox{on}~[0,1),  &\mbox{and}& u(1,x) = \f(x).
\eea
Here, $\pa_t$ and $D^2$ denote, respectively, the partial
derivative with respect to $t$,
and the partial Hessian with respect to the space variable $x$.
Then, the conditional $G$-expectation
of the random variable $\f(B_{1})$ at time $t$ is defined by
\beaa
\dbE^G_t \left[\f(B_{1})\right] &:=& u\left(t,B_{t}\right).
\eeaa
In particular,
the $G$-expectation of $\f(B_{1})$ is given by
\beaa
\dbE^G [\f(B_{1})] &:=& \dbE^G_0 [\f(B_{1})] = u(0,0).
\eeaa

Next consider the random variables of
the form $\xi:=\f (B_{t_1},\ldots,B_{t_{n-1}},B_{t_n})$
for some bounded
Lipschitz continuous function $\f$ on $\dbR^{d\times n}$
and $0\le t_1< \ldots <t_{n}=1$.
For $t_{i-1}\le t<t_{i}$, let
$$
\dbE^G_t \left[\xi\right ]=
\dbE^G_t \left[\f(B_{t_1},\cds, B_{t_n})\right ]
:=v_{i}(t,B_{t_1}, \cds, B_{t_{i-1}},B_t),
$$
where $\{v_i\}_{i=1,\ldots,n-1}$ is the unique,
bounded, Lipschitz viscosity solution of the following
equation,
\bea
\label{Gpde2}
-\pa_t v_i - G\left(D^2 v_i\right) &=& 0,\qquad t_{i-1}\le t< t_{i}
\qquad \mbox{and}\\
v_i\left(t_{i},x_1,\cds,x_{i-1},x\right) &= &v_{i+1}
\left(t_{i},x_1,\cds,x_{i-1},x,x\right),
\nonumber
\eea
and $v_n$ solves the above equation with final data
$v_n(1,x_1,\ldots, x_{n-1},x)=\f(x_1,\ldots,x_{n-1},x)$. 
Here, for $v_i$, the variables $(x_1,\cds, x_{i-1})$ are 
(fixed) parameters and the Hessian $D^2$ is the second 
order derivative on $x$. Moreover,
if we set  $u_i(x_1,\ldots,x_i)=v_{i+1}(t_i,x_1,\ldots,x_i,x_i)$, then
for $t_{i-1}\le t<t_{i}$ we have the following additional identity,
$$
\dbE^G_t \left[\f(B_{t_1},\cds, B_{t_n})\right ]
=v_{i}(t,B_{t_1}, \cds, B_{t_{i-1}},B_t)
=\dbE^G_t \left[u_{i}(B_{t_1},\cds, B_{t_{i}})\right ].
$$

Let $\cL_{ip}$ denote the space of all
random variables of the form $\f(B_{t_1},\cds, B_{t_n})$ with
a bounded and Lipschitz function $\f$. For $p\ge 1$,
$\cL^p_G$ is the closure of  $\cL_{ip}$  under the norm
\beaa
\|\xi\|_{\cL^p_G}^p := \dbE^G[|\xi|^p].
\eeaa
We may then extend the definitions of the $G$-expectation and
the conditional $G$-expectation to all $\xi\in \cL^1_G$.
In particular, the important tower property of the conditional
expectation still holds, 
\bea
\label{condG}
\dbE^G\Big[\dbE^G_t[\xi]\Big] = \dbE^G[\xi] &\mbox{for all}& \xi \in \cL^1_G.
\eea
A characterization of this space, in particular a Lusin type theorem,
is obtained in \cite{DHP}.  However, since these integrability classes
are defined through the closure of a rather smooth space $\cL_{ip}$,
they require substantial ``smoothness".  Indeed,
in the Appendix, we construct a bounded random
variable which is not in $\cL_G^1$ (see Example \ref{e.l1}).

We now can define $G$-martingales.
\begin{defn}
\label{Gmartingale}  An $\dbF$-progressively measurable $\cL^1_G$-valued process
$M$ is called a {\it{$G$-martingale}}  if and only if for any $0\le s < t$, $M_s = \dbE^G_s[M_t]$.

$M$ is called a{\it{ symmetric $G$-martingale}},
if both $M$ and $-M$ are $G$-martingales.
\end{defn}

A $G$-stochastic integral (as will be defined in the next subsection)
is an example of a symmetric $G$-martingale. In particular, the 
canonical process $B$ is a symmetric $G$-martingale.  But not all
$G$-martingales are stochastic integrals and not all are symmetric.

\subsection{Stochastic integral and quadratic variation}
\label{ss.Gintegral}

For $p \in [1,\infty)$, we let $\cH^{p,0}_G$ be the space of $\dbF$-progressively 
measurable, $\dbR^d$-valued
piecewise constant processes $H= \sum_{i\ge 0} H_{t_i}\1_{[t_i, t_{i+1})}$ such that
 $H_{t_i}\in \cL^p_G$. For $H\in \cH^{p,0}_G$, the $G$-stochastic integral is easily defined by
\beaa
\int_0^t H_s d_G B_s := \sum_{i\ge 0} H_{t_i}[B_{t\wedge t_{i+1}}-B_{t\wedge t_i}].
\eeaa
Notice that this definition is completely universal in the sense that it
 is pointwise and independent of $G$. Let $\cH^p_G$ be the 
 closure of $\cH^{p,0}_G$ under the norm:
 $$
 \|H\|_{\cH^p_G}^p := \int_0^1 \dbE^G[|H_t|^p] dt.
$$
By a closure argument the stochastic integral is defined
for all $H \in \cH^p_G$.

It is clear that the set of $G$-martingales does not form 
a linear space (unless $\underline{a} = \overline{a}$).
However, for any $H\in\cH^{p,0}_G$, one may directly verifies that
the stochastic integral process $M:=\int_0^\cd H_s d_G B_s$ is
a $G$-martingale and so is $-M$. Hence, any $G$-stochastic integral
is a symmetric $G$-martingale.

This notion of the stochastic integral can be used to define
the quadratic variation process $\la B\ra^G_t$ as well.  Indeed, the
$\dbS_d$-valued process is
defined by the identity
\bea
\label{e.gqv}
\la B\ra^G_t &:=& \frac12 B_t\otimes B_t -\int_0^t
B_s \otimes d_GB_s,\ \ \forall
\ \ 0  \le t \le 1,
\eea
where the tensor product $\otimes$
is as in the Notations \ref{ss.notations}.
We can directly check that the integrand $B_t$ is in the integration class
$\cH^p_G$.  Therefore, $\la B\ra^G_t$ is well defined.

\section{Quasi-sure stochastic analysis of Denis \& Martini \cite{DM}}
\label{s.qs}
\setcounter{equation}{0}

Let $\dbP$ be a probability measure on $(\O,\dbF)$
so that the canonical process $B$ is a martingale.
Then, the quadratic variation process $\la B\ra_t$
of $B$ under $\dbP$ exists. We consider the
subset $\cP:=\overline\cP^W_{[\underline a, \overline a]}$ of such
measures $\dbP$ so that $\la B\ra_t$ satisfies
the following for some deterministic constant $c=c(\dbP) >0$,
\bea
\label{abound}
0< \left[c I_d \vee \underline{a}\right]
\le \frac{d\la B\ra_t}{dt} \le \overline{a}, \qquad
\forall \ t\in[0,1],\  \dbP-\mbox{a.s.},
\eea
where $I_d$ is the identity matrix in $\dbS_d$.
Notice that when $\underline{a}$ is positive definite, as required in 
Denis and Martini \cite{DM},
we do not need $c I_d$ in the lower bound.
Also, the constant $c=c(\dbP)$ may be different for each measure.
Denis and Martini \cite{DM} define the following.

\begin{defn}\label{def-qs}
We say that a property holds $\cP-$quasi-surely,
abbreviated as q.s., if it holds $\dbP$-almost surely for all $\dbP\in\cP$.
\end{defn}

\begin{rem}
\label{r.strongP}
{\rm{All the results in this paper will also hold true if we let
$\cP:=\overline\cP^S_{[\underline a, \overline a]}$ be the set
of all probability measures $\dbP^\a$
given by
 \beaa
 \dbP^\a := \dbP_0 \circ (X^\a)^{-1} &\mbox{where}& X^\a_t :=
 \int_0^t \a_s^{1\slash 2} dB_s, t\in [0,1], \dbP_0-\mbox{a.s.}
\eeaa
for some $\dbF-$progressively measurable process $\a$ taking 
values in $\dbS_d$ and satisfying
$$
\left[c(\a) I_d \vee \underline{a}\right]
\le \a_t \le \overline{a}, \qquad
\forall t\in[0,1],\  \dbP_0-\mbox{a.s.},
$$
where the constant $c(\a)>0$ may depend on $\a$.
We note that $\overline\cP^S_{[\underline a, \overline a]}$ is a
strict subset of $\overline\cP^W_{[\underline a, \overline a]}$ and
each $\dbP\in \overline\cP^S_{[\underline a, \overline a]}$ satisfies
the Blumenthal zero-one law and the martingale representation property.
We remark that Denis and Martini \cite{DM} uses the space
$\overline\cP^W_{[\underline a, \overline a]}$. 
But Denis, Hu and Peng \cite{DHP}
and our subsequent work \cite{STZ09b} essentially use
$\overline\cP^S_{[\underline a, \overline a]}$.}}\ep
\end{rem}

The following are immediate consequences of the definition of
$G$-expectations.

\begin{prop}
\label{p.agree}
Let $H \in \cH^2_G$.  Then, $H$ is It\^o-integrable
for every $\dbP \in \cP$.  Moreover,
\bea
\label{integral}
\int H_sd_GB_s =
\int H_sd B_s, \ \ \ \dbP\mbox{-a.s. for every}~\dbP\in\cP,
\eea
where the right hand side is the usual It\^o integral.  Consequently,
the quadratic variation process $\la B\ra^G$ defined in \reff{e.gqv}
agrees with the usual
quadratic variation process quasi surely.
\end{prop}
\proof
The above statements clearly hold for the integrands $H\in
\cH_G^{2,0}$ (i.e. the piece-wise constant processes).
For $H\in \cH^2_G$, there exist $H^n\in \cH_G^{2,0}$ such that 
$\lim_{n\to\infty}\|H^n-H\|_{\cH^2_G}=0$. 
For any fixed $\dbP \in \cP$, since 
$\dbE^\dbP[\int_0^1|H^n_t-H_t|^2dt] \le \|H^n-H\|^2_{\cH^2_G}$,
the equality \reff{integral} holds. The statement about the quadratic variation
follows from the general statement about the stochastic integrals
and the formula \reff{e.gqv}.
\ep

Next we recall a dual characterization of the $G$-expectation
as proved in \cite{DHP}.  We will then generalize that characterization
to the $G$-conditional expectations.  Like the previous result,
this generalization is also an immediate consequence of the
previous results.  We need the following notation,
for $t \in [0,1]$ and $\dbP \in \cP$,
\bea
\label{e.ptp}
\cP(t,\dbP) := \left\{\dbP'\in \cP: ~\dbP' = \dbP ~\mbox{on}~ \cF_t\right\}.
\eea
Notice that for any $\dbP^\prime \in \cP(t,\dbP)$ and $\xi \in \cL^1_G$,
the random variable $\dbE^{\dbP^\prime}\left[ \xi  | \cF_t\right]$
is defined both $\dbP$ and $\dbP^\prime$ almost surely.
Also recall that  $\esup=\esup^{\dbP} $ is the essential supremum of
a class of $\dbP$ almost surely defined random variables.  Clearly,
it is also defined $\dbP$ almost surely
(see Definition A.1 on page 323 in \cite{KS}).  In particular,
for $t \in [0,1]$, we may define
\begin{equation}\label{extension Gexp}
{\esup_{\dbP^\prime
\in \cP(t,\dbP)}}\ \dbE^{\dbP^\prime}\left[ \ \xi \ |\ \cF_t\right]
\end{equation}
as a $\dbP$-almost sure random variable. We remark that, for given $\dbP$, 
the above random variable can be first defined as $\cF^\dbP_t$-measurable. 
However, in view of Lemma \ref{lem-filtration}, 
we will always consider its $\cF_t$-measurable version.

We now have the following characterization of the $G$-conditional
expectation.

\begin{prop}
\label{p.gcond}
For any $\xi \in \cL^1_G$, $t \in [0,1]$, and $\dbP\in\cP$,
$$
\dbE^G_t\left[ \xi \right] = \esup_{\dbP^\prime
\in \cP(t,\dbP)} \ \dbE^{\dbP^\prime}\left[ \ \xi \ |\ \cF_t\right],
\qquad \dbP-a.s..
$$
Moreover, an $\dbF$-progressively measurable $\cL^1_G$ valued process
$M$ is a $G$-martingale
if and only if it satisfies the following
dynamic programming principle for all $0 \le s \le t \le 1$
and $\dbP \in \cP$,
\bea
\label{e.dpp}
M_s = \esup_{\dbP^\prime
\in \cP(s,\dbP)} \ \dbE^{\dbP^\prime}\left[\ M_t \ |\ \cF_s\right],
\qquad \dbP-a.s..
\eea
\end{prop}
\proof
The characterization of the conditional expectation
follows  from
\cite{DHP} for $\xi \in \cL_{ip}$.
Indeed, \cite{DHP} proves this result when
the set of probability measures is $\overline\cP^S_{[\underline a, \overline a]}$
as defined in Remark \ref{r.strongP}.
Moreover when $\xi= g(B_1)$, we can use the 
dynamic programming equation \reff{Gpde1}
and classical verification arguments as in \cite{FS}
to conclude the claimed representation
in our formulation.  Then, a simple induction argument
extends the result to all $\xi \in \cL_{ip}$.

For $\xi\in \cL^1_G$,  
there exist $\xi_n\in \cL_{ip}$ such that 
$\lim_{n\to\infty}\dbE^G\left[\left|\xi_n-\xi \right|\right]=0$. 
Then, for every $t\in [0,1]$, by the definition of $\dbE^G_t[\xi]$, 
$$
\lim_{n\to\infty}\dbE^G\left[\left|
\dbE^G_t[\xi_n]-\dbE^G_t[\xi]\right|\right]=0.
$$ 
Moreover, for any $t\in [0,1]$ and $\dbP\in \cP$, 
$$
\dbE^\dbP\left[\left|\dbE^G_t[\xi_n]-\dbE^G_t[\xi]\right|\right]
\le \dbE^G\left[\left|\dbE^G_t[\xi_n]-\dbE^G_t[\xi] \right|\right].
$$
Using these and \reff{condG},
we directly estimate that
\beaa
\dbE^\dbP\left[\left|\esup_{\dbP^\prime
\in \cP(t,\dbP)} \ \dbE^{\dbP^\prime}_t[\xi_n]-\esup_{\dbP^\prime
\in \cP(t,\dbP)} \ \dbE^{\dbP^\prime}_t[\xi]\right|\right]
&\le& \dbE^\dbP\Big[\esup_{\dbP^\prime
\in \cP(t,\dbP)} \ \dbE^{\dbP^\prime}_t[|\xi_n-\xi|]\Big]\\
&\le& \dbE^\dbP\Big[\esup_{\dbP^\prime
\in \cP(t,\dbP)} \ \dbE^G_t[|\xi_n-\xi|]\Big]= \dbE^\dbP\Big[\dbE^G_t[|\xi_n-\xi|]\Big]\\
&\le& \dbE^G\Big[\dbE^G_t[|\xi_n-\xi|]\Big] = \dbE^G[|\xi_n-\xi|].
\eeaa
Therefore,
\beaa
\dbE^G_t[\xi] = \lim_{n\to\infty} \dbE^G_t[\xi_n] = \lim_{n\to\infty} \esup_{\dbP^\prime
\in \cP(t,\dbP)} \ \dbE^{\dbP^\prime}_t[\xi_n] =  \esup_{\dbP^\prime
\in \cP(t,\dbP)} \ \dbE^{\dbP^\prime}_t[\xi],~~\dbP\mbox{-a.s..}
\eeaa
The martingale property
is a direct consequence of the tower property
of the $G$-conditional expectation as proved in
\cite{Peng1} and the above formula for the conditional expectation.
\ep

\begin{rem}
\label{r.nicole}
{\rm{
In their classical paper \cite{EJ}, El Karoui \& Jeanblanc consider a very general
stochastic optimal control problem.  Their results in our context imply that
$$
M_s^\dbP := \esup_{\dbP^\prime
\in \cP(s,\dbP)} \ \dbE^{\dbP^\prime}\left[\xi \ |\ \cF_s\right]
$$
is a $\dbP$-super-martingale for all $\dbP \in \cP$.  Moreover $\dbP^*$ is
a maximizer if and only if $M^{\dbP^*}$ is
$\dbP^*$-martingale.  While this result provides a characterization
of the optimal measure $\dbP^*$, it does not provide a ``universal"
hedge.  More precisely their approach
provides an optimal control which is defined
only for the optimal measure and on its support.
Indeed, the super-martingale property of $M^\dbP$ imply that
there are an increasing process $K^\dbP$ and an integrand $H^\dbP$
so that
$$
M^\dbP_t= \int_0^t H^\dbP_s dB_s - K^\dbP_t.
$$
However, aggregating these processes into one universally defined
$K$ and $H$ is not immediate.  In the standard Markovian
context, this problem can be solved directly.  However, it is
exactly the non-Markovian generalization that motivates
this paper and \cite{DM,Peng,Peng1}.
This interesting question of aggregation is further
discussed in the Remark \ref{r.aggregation}.}}\ep
\end{rem}

\section{Spaces and Norms}
\label{s.norms}
\setcounter{equation}{0}

The particular case of $t=0$ in \reff{e.dpp} gives the following
dual characterization proved in \cite{DHP},
$$
\dbE^G \left[ \xi \right] = \sup_{\dbP\in \cP}
\ \dbE^{\dbP}\left[  \xi  \right].
$$
The above results enable us to extend the definition of
$G$-expectation and $G$-martingales
to a possibly larger class of random variables.  In particular,
this extension has the advantage of not referring to the
partial differential equation \reff{Gpde1}.
We will not
develop this theory here.
However, in view of the results and the norms used in the theory
of BSDEs, we introduce the following  function
spaces.

For $p\ge 1$, and an $\cF_1$-measurable, non-negative random
variable $\xi$, we set
$$
\|\xi \|_{\dbL^p_\cP}^p := \sup_{\dbP \in \cP} \
\dbE^\dbP \left[ \esup_{t \in [0,1]} \ \left( M_t^\dbP(\xi) \right)^p \right],\ \
{\mbox{where}}\ \
M_t^\dbP(\xi) := \esup_{\dbP^\prime \in \cP(t,\dbP)}
\ \dbE^{\dbP^\prime} \left[ \xi |\cF_t\right].
$$

In the above definition,  {\em{a priori}} we do not have any 
information on the regularity of $M_t^\dbP(\xi)$ on its $t$
dependence.   That is the reason for
defining the norm through the random variable
$ \esup_{t \in [0,1]} \ \left( M_t^\dbP(\xi) \right)$,
which is,  in view of Lemma \ref{lem-filtration}, 
$\cF_1$-measurable.  Alternatively, one may first prove that
$M_t^\dbP(\xi)$ is a $\dbP$-supermartingale and that it admits a {\cad} version.
Then, $\sup_{t \in [0,1]}\ \left( M_t^\dbP(\xi) \right)^p$ 
would be  measurable and we could use it in the definition.
However, we believe that this issue tangential to the main 
thrust of the paper and we prefer to give the above 
quicker definition.

We next define
\bea
\label{e.lnorm}
\dbL^p_\cP &: = &\left\{ \xi  \ :\  \cF_1\mbox{-measurable and}~ \| \xi  \|_{\dbL^p_\cP}:=\| |\xi|  \|_{\dbL^p_\cP} < \infty \right\},\\
\nonumber
\cL^p_\cP&:=& {\mbox{closure of}}\ \cL_{ip}\ {\mbox{under the  norm}}\  \dbL^p_\cP.
\eea
Notice that if $\xi \in \cL^1_G$, then
$M_t^\dbP(\xi)= \dbE^G_t[\xi]$
for every $\dbP \in \cP$.  Moreover, for every $\xi\in \cL_{ip}$,
$\|\xi \|_{\dbL^p_\cP}=\|\xi^*\|_{\cL^p_G}$, where $\xi^*:=\sup_{t \in [0,1]}\  \dbE^G_t \left[  |\xi |\right] $.

In the Appendix, we compare the integrability
classes defined by Peng \cite{Peng} and the above spaces.
The connection is related to the Doob maximal inequalities in the
setting of  $G$-expectations.  In particular, we prove the following.
\begin{lem}
\label{l.equal}
${\cup_{p>2}} \cL^p_G \subset \cL^2_\cP \subset \dbL^2_\cP \cap \cL^2_G \subset \dbL^2_\cP$.
Moreover, the final inclusion is strict.
\end{lem}

We also define the following norms for the processes.
As usual $1 \le p<\infty$.  For an $\dbF$-progressively measurable integrand $H$ and
a stochastic process $Y$,
we set
\bea
\label{e.hnorm}
\|H\|^p_{\dbH^p_\cP}&:=&  \sup_{\dbP\in \cP} \dbE^{\dbP}\left[\left
(\int_0^1 (d\la B\ra_tH_t\cdot H_t)\right)^{p\over 2}\right],\\
\|Y\|^p_{\dbS^p_\cP}& := & \sup_{\dbP\in\cP} \dbE^\dbP\Big[\esup_{0\le t\le 1} |Y_t|^p\Big].
\eea
If $Y_t= \dbE^G_t[|\xi |]$ for some $\xi \in \cL^1_G$, then $\|Y\|^p_{\dbS^p_\cP} =
\|\xi \|_{\dbL^p_\cP}^p$.  This identity also motivates the definition of the  norm
$\dbL^p_\cP$.  Moreover, when the lower bound $\underline{a}$ in \reff{abound} is non-degenerate,
then the $\dbH^p_\cP$ norm is equivalent to the norm used in \cite{DHP,Peng1}:
$$
\sup_{\dbP\in \cP} \dbE^{\dbP}
\left[\left(\int_0^1 |H_t|^2dt\right)^{p\over 2}\right].
$$
In analogy with the standard notation in stochastic calculus, we define the following
spaces.
\begin{defn}
\label{d.hnorm}{\rm{Let $p \in [1,\infty)$ and $\cP$ be
as in Section \ref{s.qs}.}}
\begin{itemize}
\item {\rm{ $\dbH^p_{\cP}$ is the set of all $\dbF$-progressively measurable integrands with a finite
$\| \cdot \|_{\dbH^p_\cP}$-norm,}}
\item {\rm{$\cH^p_\cP$ is the closure of ${\cH^{p,0}_G}$
under the norm $\| \cdot \|_{\dbH^p_\cP}$,}}
\item {\rm{$\dbS^p_\cP$  is the set of all $\dbF$-progressively measurable  processes
 with quasi surely {\em{ continuous}} paths and
finite $\|\cdot\|_{\dbS^p_\cP}$-norm,}}
\item {\rm{ $\dbI^p_\cP$ is the subset of $\dbS^p_\cP$ of  {\em{non-decreasing}} processes with
$X_0=0$.}}
\ep
\end{itemize}
\end{defn}
Clearly all of the above spaces are defined as quasi-sure equivalence classes.
As such, they are complete and therefore Banach spaces.
Also $\|H\|_{\dbH^p_\cP} \le \|H\|_{\cH^p_G}$ for $H\in \cH^{p,0}_G$,
then it is clear that $\cH^p_G \subset \cH^p_\cP$.
Therefore $\cH^p_\cP$ is the
 closure of $\cH^p_G$ under the norm $\|\cd\|_{\dbH^p_\cP}$.

\begin{rem}
\label{r.aggregation}
{\rm Given $\xi\in \dbL^1_\cP$ (but not necessarily in $\cL_G^1$) and an
$\dbF$-stopping time $\t$,
it is not straightforward to define the conditional $G_\cP-$expectation
$\dbE^\cP_{\t}[\xi]$ as in \reff{extension Gexp}. Indeed, set
\beaa
M^\dbP_{\t}
&:=& \esup_{\dbP' \in \cP(\t, \dbP)} \dbE^{\dbP'}_{\t}[\xi],
~~\dbP-\mbox{a.s.}
\eeaa
Then, to define the conditional expectation, we need to aggregate this family
of random variables  $\{M^\dbP_{\t}, \dbP\in\cP\}$ into one ``universally" defined
random variable.  A similar problem arises in the definition of a stochastic integral
for a given integrand  $H\in  \dbH^2_{\cP}$.  Again, for $\dbP \in \cP$, we set
$M^\dbP_t:= \int_0^t H_s dB_s$.  Then, to define the $G$-stochastic integral of
$H$ we need to aggregate this family
of stochastic processes.

The issue of aggregation is an interesting
technical question.  Generally, a solution to this technical issue
is given by imposing regularity on the random variables.
Indeed, for all random variables
which are in $\cL^p_G$, one can define the
universal version through a closure argument.
However, there are other alternatives and a
comprehensive study of this question is given
in our accompanying paper  \cite{STZ09a}.

Finally we recall that, when the integrand $H$ has the additional
regularity that it is a  {\cad}  process, then
Karandikar \cite{kar} defines the stochastic integral
$M^\dbP_t:= \int_0^t H_s dB_s$  point wise.  This definition
can then be used as the aggregating process.   \ep
} \end{rem}

\section{The martingale representation theorem}
\label{s.mrt}
\setcounter{equation}{0}

To motivate the main result of this paper, we
first consider the case $\xi=\varphi(B_1)$ for some
smooth, bounded function $\varphi$. In this case, as in Peng \cite{Peng-mono,Peng1},
a formal construction can be derived by simply using the It\^o's formula.
Now suppose that the solution $u(t,x)$ of \reff{Gpde1} is smooth.
Indeed, we can approximate the equation
\reff{Gpde1} so that the approximating equation
admits smooth solutions as proved by \cite{krylov}.  This is done in the
Appendix.
Then, we set $Y_t:= u(t,B_t)=\dbE^G_t[\xi] $, $H_t:= \nabla u(t,B_t)$ and
$$
K_t := \int_0^t \left(G(D^2u(s,B_s)) - \frac12 \tr
\left[\hat a_sD^2u(s,B_s)\right]\right)ds,\qquad
\hat a_t := \frac{d\la B \ra_t }{dt}, \ q.s..
$$
Using \reff{Gpde1}, \reff{abound} and
the definition of the nonlinearity $G$,
one may directly check that
$$
Y_t = \xi - \int_t^1 H_s dB_s + K_1-K_t,\ \
{\mbox{and}}\ \ dK_t \ge 0 \qquad q.s..
$$
Also, the characterization of $G$-martingales in Proposition \ref{p.gcond}
and the definition of the nonlinearity $G$ imply that $-K$ is a $G$-martingale.
Hence for the random variable $\xi=\varphi(B_1)$,
we have the martingale representation.
More importantly, this example also shows that in general
a non-decreasing process $K$ is always present in this representation.
The above construction is also the basic step in our construction.
Indeed essentially for almost all
random variables in $\cL_{ip}$ the above construction
proves the result.  We then prove that stochastic integrals and
non-decreasing martingales are closed subsets under the appropriate
norms as defined in the preceding section.
Finally, these results  allow us to prove the result by a closure argument.

\subsection{Main results}
\label{ss.main}

We first state the main result.  Recall
that function spaces are defined in Definition \ref{d.hnorm}.
\begin{thm}
\label{t.main}  Assume that  $\underline{a}$ and $\overline{a}$
satisfy \reff{e.matrix}.  Then,
for every $\xi \in \cL^2_\cP$, the conditional $G-$expectation
process $Y_t := \dbE^G_t[\xi]$ is in $\dbS^2_\cP$,
and there exist unique  $H \in \cH^2_\cP$,
$K \in \dbI^2_\cP$ so that $N:=-K$ is a $G$-martingale
and for every $t \in [0,T]$,
\be
\label{e.main}
Y_t  = \xi -\int_t^1 \ H_s dB_s +K_1 - K_t\
= \dbE^G[\xi] + \int_0^t \ H_s dB_s - K_t,\qquad
q.s..
\ee
In particular, the stochastic integrals are defined both as $G$-stochastic
integrals and also quasi surely.  Moreover the following
estimate is also satisfied with a universal constant $C^*$,
\be
\label{e.est}
\|Y\|_{\dbS^2_\cP} + \|H\|_{\dbH^2_\cP} + \|K\|_{\dbS^2_\cP} \le
C^* \|\xi\|_{\dbL^2_\cP}.
\ee
\end{thm}

The proof of the above theorem will be completed in several lemmas below.

In the above theorem the integrand $H$ is not only in the
class $\dbH^2_\cP$ but also in the closure of $\cH^{2,0}_G$
under the norm $\|\cdot \|_{\dbH^2_\cP}$.  Indeed this fact implies that
stochastic integral is well defined quasi surely as it is shown in the
next subsection.

The following is an immediate corollary of the above martingale representation.

\begin{cor}
\label{c.symmetric}
A $G$-martingale $M$ with $M_1\in \cL^2_\cP$ is symmetric if and only
if the process $K$ in the representation \reff{e.main}
is identically equal to zero.
\end{cor}

In addition to the estimate \reff{e.est} an estimate of the differences of the solutions
is known to be an important tool. Let $\xi_1, \xi_2\in \cL^2_{\cP}$ and $(Y^i, H^i, K^i)$ be
the processes in the martingale representation. We set
$\d \xi:= \xi^1-\xi^2$, $\d Y:= Y^1-Y^2$,
$\d Z:= Z^1-Z^2$ and $\d K:= K^1-K^2$.

\begin{thm}
\label{t.difference}
There exists a  universal constant $C^*$ so that,
\beaa
\|\d Y\|_{\dbS^2_\cP} &\le& \|\d \xi\|_{\dbL^2_{\cP}},\\
\|\d H\|_{\dbH^2_\cP} + \|\d K\|_{\dbS^2_\cP}
&\le& C^* \left [\|\d \xi\|_{\dbL^2_{\cP}}+
\left(\|\xi^1\|_{\dbL^2_\cP}^{1\over 2} + \|\xi^2\|_{\dbL^2_\cP}^{1\over 2}
\right) \ \|\d \xi\|_{\dbL^2_{\cP}}^{1\over 2}\right].
\eeaa
\end{thm}

\subsection{Stochastic Integral and Symmetric G-martingales}
\label{ss.symmetric}

As discussed in Remark \ref{r.aggregation},
for an integrand $H\in \dbH^2_\cP$ it is not immediate to
define the stochastic integral $\int_0^\cdot H_s dB_s$ quasi surely.
However, the stochastic integral is defined in \cite{Peng} for integrands
$H \in \cH^{2,0}_G$.
Then,
for integrands in $\cH^2_\cP$
a closure argument can be used to construct the stochastic integral quasi-surely.
(Recall that $\cH^2_\cP$ is the closure of ${\cH^{2,0}_G}$
under the norm $\| \cdot \|_{\dbH^2_\cP}$.  )

\begin{thm}
\label{t.integral}
For any $H\in \cH^2_\cP$, the stochastic integral  $\int_0^\cdot H_s dB_s$
exists quasi surely.  Moreover, the
stochastic integral satisfies the Burkholder-Davis-Gundy inequality
\be
\label{e.bdg}
\|H\|_{\dbH^2_\cP}\le
\Big\|\int_0^\cd H_s dB_s\Big\|_{\dbS^2_\cP}
\le 2 \|H\|_{\dbH^2_\cP}.
\ee
\end{thm}

\proof Let $H\in \cH^2_\cP$.  Then, there is a sequence
$\{H^n\}_{n} \subset \cH^{2,0}_G$ so that
$\| H^n-H\|_{\dbH^2_\cP}$ converges to zero as $n$ tends to infinity.
By relabeling the sequence we may assume that
$ \|H^n-H\|_{\dbH^2_{\cP}} \le 2^{-n}$ for every $n$.
Moreover, since $ H \in\dbH^2_{\cP}$, for every $\dbP \in \cP$,
$$
M^\dbP_t := \int_0^t H_s dB_s, \qquad t\in[0,1] ,
$$
is $\dbP$-almost surely well-defined.
Since $H^n \in \cH^{2,0}_G$, the $G$-stochastic integral
$$
M^n_t := \int_0^t H^n_s dB_s,~~t\in[0,1],
$$
is also  defined pointwise.

We now have to prove that the family
$\{M^\dbP,P\in\cP\}$ can be aggregated into a universal 
$\dbF$-progressively measurable  process. For this,
we define
$$
\overline M_t := \limsup_{n\to\infty} M^n_t,~~t\in[0,1].
$$
Notice that $\overline M$ is pointwisely defined and 
$\dbF$-progressively measurable.
We continue by showing that $\overline M = M^\dbP $,
$\dbP$-almost surely, for every $\dbP \in \cP$.
Indeed for any $\dbP\in\cP$,
we use the Burkholder-Davis-Gundy
inequality to obtain
\beaa
\dbE^\dbP\Big[\sup_{0\le t\le 1} |M^n_t - M^\dbP_t|^2\Big]
& = &
\dbE^\dbP\Big[\sup_{0\le t\le 1} |\int_0^t (H^n_s - H_s)dB_s|^2\Big]\\
&\le&
4 \dbE^\dbP\Big[ |\int_0^1 (H^n_s - H_s)dB_s|^2\Big]\\
&=&
4 \dbE^\dbP\Big[\int_0^1 |\hat a_s^{1\slash 2}(H^n_s - H_s)|^2ds\Big]\\
& \le &  4 \|H^n-H\|_{\dbH^2_{\cP}}^2  \le 2^{2-2n}.
\eeaa
We then directly estimate that
\beaa
\sum_{n=1}^\infty \dbP\Big[\sup_{0\le t\le 1} |M^n_t - M^\dbP_t| 
\ge n^{-2}\Big] &\le&
\sum_{n=1}^\infty n^2\dbE^\dbP\Big[\sup_{0\le t\le 1} 
|M^n_t - M^\dbP_t|^2\Big]^{1\over 2} \;<\;
\infty.
\eeaa
By the Borel-Cantelli Lemma,
$$
\lim_{n \to \infty}  \sup_{0\le t\le T} |M^n_t - M^\dbP_t|=
0,~~\dbP-\mbox{a.s.}\ .
$$
This implies that $M^\dbP = \overline M$, $\dbP-$almost surely.
Since this holds for every $\dbP\in\cP$, we conclude that
the process $\overline M$ is an aggregating process.  Hence the
stochastic integral is defined.

The Burkholder-Davis-Gundy inequalities follow directly from the definitions.
\ep

\vspace{5mm}

We close this subsection by stating the
following result for symmetric $G$-martingales, which is
 an immediate consequence of
the main results.
\begin{thm}
\label{thm-symmetric} Let $M$ be a $G$-martingale with $M_1\in \cL^2_\cP$. 
The following are equivalent:\\
{\rm (i)} $M$ is a $\dbP-$martingale for every $\dbP\in\cP$,
\\
{\rm (ii)} $M$ is a symmetric $G-$martingale,
\\
{\rm (iii)} For any $G-$martingale $N$, both $N+M$ and $N-M$
 are also $G-$martingales,
\\
{\rm (iv)} $\dbE^G\{-M_t\} = - \dbE^G\{M_t\}$ for any $t\ge0$,
\\
{\rm (v)} There exists
 $H\in \cH^2_\cP$ so that $M_t := M_0+\int_0^t H_s dB_s$.
\end{thm}

\begin{rem}
\label{rem-main}
{\rm{The main reason for the requirements 
$\xi\in\cL^2_\cP$ and $H\in \cH^2_\cP$ is to ensure 
the existence of the universal version of the conditional 
$G$-expectation $E^G_t[\xi]$ and the stochastic integral $\int_0^t H_s dB_s$. 
However, if we are given a $G$-martingale $M$ with 
$M_1 \in \dbL^2_\cP$, then there would be  no aggregation issue. 
Then, following the same arguments, one can easily show that  
Theorem \ref{thm-symmetric} still holds true
under the weaker assumption $M_1\in \dbL^2_\cP$.
Moreover, (v) requires only $H\in \dbH^2_\cP$.}}
\end{rem}

Recall that $ \dbI^2_\cP$ is defined in Definition \ref{d.hnorm} as
the set of all $\dbF$-progressively measurable, non-decreasing, 
continuous processes with
finite $\|\cdot\|_{\dbS^p_\cP}$.  For $(H, K) \in \cH^2_{\cP}\times \dbI^2_\cP$,
define  a process by
\be
\label{e.rep}
M_t := M_0 + \int_0^t H_s dB_s - K_t.
\ee
An immediate corollary of the above result is
the following.

\begin{cor}
\label{c.characterization}
The process $M$ defined in \reff{e.rep} is a $G$-martingale
if and only if the non-increasing process $-K$ is a $G$-martingale.
\end{cor}

\subsection{Increasing $G$-martingales}
\label{ss.increasing}

In this section we show that the set of non-decreasing
$G$-martingales is a closed set.
Indeed, let ${MI}^2_{\cP}$ be the set of all processes
$K\in \dbI^2_{\cP}$ such that $-K$ is a $G$-martingale.
Then we have the following closure result which is similar to
Theorem \ref{t.integral}.
\begin{thm}
\label{t.Kclose}
The space $MI^2_{\cP}$ is closed in $\dbS^2_{\cP}$
under norm $\|\cd\|_{\dbS^2_{\cP}}$.
\end{thm}

\proof
Consider a sequence $K^n \in MI^2_{\cP}$ converging
to a process $K\in\dbI^2_{\cP}$ in the norm $\|\cd\|_{\dbS^2_{\cP}}$.
We claim that the limit $-K$ is also a $G$-martingale
and therefore $K \in   MI^2_{\cP}$.
Indeed, for every $0 \le s\le t\le 1$, set $A_t:=K_t-K_s$ 
and $A^n_t:=K^n_t-K^n_s$.
Then, by the martingale property of the sequence,
for every $n$ and $\dbP \in \cP$, we have
\beaa
\einf_{\dbP'\in\cP(s,\dbP)}\dbE^{\dbP'}_s[A^n_t] =0,
&\dbP-\mbox{a.s.}.
\eeaa
Moreover, $\dbP-$a.s.,
$$
\einf_{\dbP'\in\cP(s,\dbP)}\dbE^{\dbP'}_s[A_t]
\le
\esup_{\dbP'\in\cP(s,\dbP)}\dbE^{\dbP'}_s|A_t-A^n_t|
+\einf_{\dbP'\in\cP(s,\dbP)}\dbE^{\dbP'}_s[A^n_t]
=
\esup_{\dbP'\in\cP(s,\dbP)}\dbE^{\dbP'}_s|A_t-A^n_t|.
$$
The following can be shown directly from the definitions:
$$
\sup_{\dbP\in\cP} \dbE^\dbP
\left[\esup_{\dbP'\in\cP(s,\dbP)}\dbE^{\dbP'}_s|A_t-A^n_t|\right]
\le
\|A-A^n\|_{\cS^2_{\cP}}.
$$
Hence by the convergence of  $\|A-A^n\|_{\cS^2_{\cP}}$ 
 to zero as $n$ tends to infinity,
we conclude that
$$
\lim_{n\to \infty} \esup_{\dbP'\in\cP(s,\dbP)}
\dbE^{\dbP'}_s|A_t-A^n_t|= 0,
\qquad \dbP-a.s.\ .
$$
Since $0\le s\le t\le 1$ and $\dbP\in\cP$ are arbitrary,
the limit process $-K$ is also a $G-$martingale.
\ep

\subsection{Estimates}
\label{ss.estimates}
For $(H, K) \in \dbH^2_{\cP}\times \dbI^2_\cP$,  let $M$ be
defined as in \reff{e.rep}.
In this subsection,
we prove certain estimates for $H$ and $K$ in terms
of the process $M$.
These  estimates are similar to those obtained
for reflected backward stochastic differential equations
in \cite{EKPPQ}.

\begin{prop}
\label{p.est1}
Let $H,K,M$ be as in \reff{e.rep}.
There exists a constant $C$ depending only on the dimension so that
$$
\|H\|_{\dbH^2_{\cP}}
+ \|K\|_{\dbS_\cP^2} \le C \|M\|_{\dbS^2_\cP}.
$$
\end{prop}

\proof
We directly calculate that
$$
d|M_t|^2 = 2  M_t H_t dB_t - 2 M_t d K_t + d\la B\ra_t H_t \cdot H_t.
$$
We integrate over $[t,1]$ to obtain,
$$
|M_t|^2 +
\int_t^1  d\la B\ra_s H_s \cdot H_s
= |M_1|^2
+ 2 \int_t^1 M_s d K_s -2 \int_t^1 M_s H_s dB_s.
$$
We then take the expected value under an arbitrary
$\dbP \in \cP$ to arrive at
$$
\dbE^\dbP\left[ |M_t|^2 +
\int_0^1  d\la B\ra_t H_t \cdot H_t \right]
\le \dbE^\dbP\left[|M_1|^2
+ 2 \int_0^1 |M_t| d K_t\right].
$$
Since $dK_t \ge 0$, for any $\e> 0$,
we have the following estimate,
\bea
\nonumber
\dbE^\dbP\left[ |M_t|^2 +
\int_0^1  d\la B\ra_t H_t \cdot H_t \right]
&\le & \dbE^\dbP\left[|M_1|^2
+ 2 \left( \sup_{t \in [0,1]} |M_t|\right) K_1\right]\\
\label{e.1}
&\le& (1+\e^{-1})  \dbE^\dbP\left[\sup_{t \in [0,1]} |M_t|^2\right]
+ \e \dbE^\dbP\left[K_1^2\right].
\eea
Next we estimate $K$.  Recall that $0=K_0 \le K_t$.
By the definition of $M_t$,
\begin{eqnarray*}
K_1^2 &=& \left(M_1-M_0
-\int_0^1 H_s dB_s \right)^2\\
&\le &
3 |M_1|^2 + 3|M_0|^2
+3 \left(\int_0^1 H_s dB_s \right)^2.
\end{eqnarray*}
We now use \reff{e.1} with $\e= {1\over 6}$.  The result is
\begin{eqnarray*}
\dbE^\dbP\left[K_1^2\right] &\le&
\dbE^\dbP\left[3 |M_1|^2 +3 |M_0|^2 +3 \int_0^1
d\la B\ra_t H_t \cdot H_t \right]\\
&\le & 27\  \dbE^\dbP\left[\sup_{t \in [0,1]} |M_t|^2\right]
+ \frac{1}{2} \dbE^\dbP\left[K_1^2\right].
\end{eqnarray*}
Hence,
$$
\dbE^\dbP\left[K_1^2\right]
\le 54\  \dbE^\dbP\left[\sup_{t \in [0,1]} |M_t|^2\right].
$$
This together with \reff{e.1} and the definitions of the norms imply the result.
\ep
\vspace{10pt}

Next we prove an estimate for differences.
So for any $(H^i,K^i)\in\dbH^2_{\cP}\times \dbI^2_{\cP}, i=1,2$,
let $M^i$ be defined as in \reff{e.rep}.  As before, let
$\d M:=M^1-M^2$, $\d H:=H^1-H^2$, $\d K:=K^1-K^2$.
\begin{prop}
\label{p.est2}
There exists a constant $C$ depending only on the dimension so that
\begin{equation}
\label{e.est2}
\|\d H\|_{\dbH^2_{\cP}}^2 + \|\d K\|^2_{\dbS^2_{\cP}}
\le C \Big[\|\d M\|^2_{\dbS^2_{\cP}}+\|\d M\|_{\dbS^2_{\cP}}
\left(\|K^1\|_{\dbS^2_{\cP}}+ \|K^2\|_{\dbS^2_{\cP}}\right)  \Big].
\end{equation}
\end{prop}
\vspace{10pt}

The terms $\| K^i\|_{\dbS^2_{\cP}} $ in the
above inequality can be estimated using Proposition \ref{p.est1}.

\proof The arguments are very similar to the proof of
Proposition \ref{p.est1}.  The only difference
is the fact that $\d K$ is no longer a monotone function.
We directly compute that
$$
\d M_t = \d M_0 + \int_0^t \d H_s dB_s - \d K_t.
$$
Then we proceed as in the proof of the previous proposition to arrive at
$$
\dbE^\dbP\left[ |\d M_t|^2 +
\int_0^1  d\la B\ra_t \d H_t \cdot \d H_t \right]
\le  \dbE^\dbP\left[|\d M_1|^2 \right]
+ \dbE^\dbP\left[ \int_0^1 |\d M_s| d |\d K|_s\right].
$$
The last integral term is directly estimated as follows.
\beaa
\dbE^\dbP\left[ \int_0^1 |\d M_s| d |\d K|_s\right]
& \le &\dbE^\dbP\left[\left(\sup_{t \in [0,1]} |\d M_t|\right)
\left(\sup_{t \in [0,1]} [|K^1_t|+|K^2_t|]\right)\right]\\
&\le& 2
\left[\dbE^\dbP\sup_{t \in [0,1]} |\d M_t|^2\right]^{1/2}
\left( \sum_{i=1}^2 \left[\dbE^\dbP\sup_{t \in [0,1]} |K^i_t|^2\right]^{1/2}\right)
\\ &\le& 2 \|\d M\|_{\dbS^2_{\cP}}
\left(\|K^1\|_{\dbS^2_{\cP}}+ \|K^2\|_{\dbS^2_{\cP}}\right).
\eeaa
The estimate of $\|\d K\|_{\dbS^2_{\cP}}$ is obtained
exactly as in the proof
of Proposition \ref{p.est1}
\ep

\subsection{Proof of Theorem \ref{t.main}}
\label{ss.proof}
We prove uniqueness first.  Suppose that there are two
pairs $(H^i,K^i)$ satisfying \reff{e.main}.
Then, we can use Proposition \ref{p.est2}
with $M^i_t=Y_t = E^G_t[\xi]$.  In particular,
$\d M \equiv 0$.  By \reff{e.est2}, we conclude that
$\|\d H\|_{\dbH^2_\cP}= \|\d K\|_{\dbS^2_\cP}=0$.

For the existence, let $\cM$ be the subset of $\dbL^2_\cP$ so that
the martingale representation \reff{e.main}
holds for all $\xi \in \cM$.
We will prove  the result by showing that
$\cM$ is closed in $\dbL^2_\cP$ and that $\cL_{ip} \subset
\cM$.  The second statement is proved in the Appendix,
by an approximation argument.  This is Proposition
\ref{p.approximate}.  Then for $\xi  \in \cL^2_\cP$
these two statements imply the existence
of $(H,K)$ as $\cL^2_\cP$
is in the closure of $\cL_{ip}$ under the norm $\dbL^2_\cP$.

To show that $\cM$ is closed, consider a sequence
$\xi^n \in \cM$ converging to $\xi \in \dbL^2_\cP$.
Since $\xi^n \in \cM$, there are $H^n \in \cH^2_\cP$
and $K^n \in \dbI^2_\cP$ so that
\reff{e.main} holds for each $n$ and $N^n:= -K^n$
is a continuous, non-increasing $G$-martingale.  We now use the
estimate \reff{e.est2} with $M^1=Y^n$ and $M^2=Y^m$
for arbitrary $n$ and $m$.
The identity
$Y^n_t=E^G_t[\xi^n]$ together with the definition of the conditional
expectation $\dbE^G_t$ imply that for every $t \in [0,1]$,
$$
\left| Y^n_t - Y^m_t\right|^2
\le \dbE^G_t\left[\left| \xi^n-\xi^m\right|^2\right].
$$
Hence the definition of the norm
$\| \cdot \|_{\dbL^2_\cP}$ yield,
$$
\| Y^n-Y^m\|_{\dbS^2_\cP}
\le \|\xi^n-\xi^m\|_{\dbL^2_\cP}.
$$
We now use the results of Propositions \ref{p.est1}
and \ref{p.est2} with $M^1=Y^n$ and $M^2=Y^m$.
The Proposition \ref{p.est1} yields for each $n$,
$$
\|K^n\|_{\dbS^2_\cP} \le \|\xi^n\|_{\dbL^2_\cP} \le
c_0:= \sup_m \|\xi^m\|_{\dbL^2_\cP}  <\infty.
$$
We use this in \reff{e.est2}. The result is
$$
\| H^n-H^m\|^2_{\dbH^2_\cP}
+ \|K^n-K^m\|^2_{\dbS^2_\cP} \le
C^*\left[ \|\xi^n-\xi^m\|^2_{\dbL^2_\cP}
+ 2c_0 \|\xi^n-\xi^m\|_{\dbL^2_\cP} \right].
$$
Hence $\{H^n\}_n$ is a Cauchy sequence
in $\cH^2_\cP$.  Therefore by the definition
of $\cH^2_\cP$, we know that
there is a limit $H\in \cH^2_\cP$.
Moreover, by \reff{e.bdg} the corresponding
stochastic integrals converge in $\dbS^2_\cP$.
Also $\{K^n\}_n$ is a Cauchy sequence
in $\dbS^2_\cP$.  By Theorem \ref{t.Kclose},
we conclude that there is a limit $K\in \dbI^2_\cP$
so that $N:=-K$ is a $G$-martingale.
Since $(Y^n,H^n,K^n)$ satisfies
\reff{e.main} with final data $Y^n_1=\xi^n$, we conclude that
the limit processes $(Y,H,K)$
also satisfies
\reff{e.main} with final data $Y_1=\xi$.
Hence $\cM$ is closed under the norm $\dbL^2_\cP$.
\ep

\subsection{Proof of Theorem \ref{t.difference}}
\label{ss.proof1}
Since $Y^i_t=E^G_t[\xi^i]$, the dual representation
of the $G$-conditional expectation yield that for each $t\in[0,1]$ ,
$$
|\d Y_t| =\left| E^G_t[\xi^1]-E^G_t[\xi^2]\right|
\le  E^G_t\left[\left| \xi^1-\xi^2\right|\right].
$$
Hence,
$$
\|\d Y\|_{\dbS^2_\cP} \le \|\d \xi\|_{\dbL^2_{\cP}}.
$$
We now use Proposition \ref{p.est2}.  The result is
$$
\|\d H\|_{\dbH^2_\cP} + \|\d K\|_{\dbS^2_\cP}
\le C^* \left[\|\d Y\|_{\dbS^2_{\cP}}+\|\d Y\|_{\dbS^2_{\cP}}^{1\over 2}\left(
\|K^1\|^{1\over 2}_{\dbS^2_\cP}+\|K^2\|^{1\over 2}_{\dbS^2_\cP}\right)\right].
$$
We now use the estimate \reff{e.est}
in the above inequality, together with the fact that $|\|\xi^2\|_{\dbS^2_\cP} - 
 \|\xi^1\|_{\dbS^2_\cP}| \le \|\d\xi\|_{\dbS^2_\cP}$, to complete
the proof of the Theorem.
\ep

\section{Appendix}
\label{s.appendix}
\setcounter{equation}{0}

In this Appendix, we construct smooth approximations of
the partial differential equations \reff{Gpde1},
\reff{Gpde2}  and study the
properties of the  integrability class $\dbL^2_\cP$.

\subsection{Approximation}
\label{ss.approximation}
The main goal of this subsection
is to construct a smooth approximation
of solutions of \reff{Gpde2}.  We require smoothness
of these solutions in order to be able to apply the It\^o rule.
The first obstacle to regularity is the possible degeneracy
of the nonlinearity $G$ or equivalently the possible degeneracy
of the lower bound $\underline a$.  Therefore,
we do not expect the equation to regularize the
final data. However, even in this case
the solution remains twice differentiable
provided the final data has this regularity.
But the second difficulty in proving smoothness
emanates from the fact that the equation \reff{Gpde2}
is solved in several time intervals and in each interval
$(t_i,t_{i+1})$ and the value $B_{t_i}$ enters into the equation as a
parameter.  Differentiability with respect to these types
of parameters is harder to prove.  Given these difficulties,
we approximate the equation as follows.

For $\epsilon \in (0,1]$, set $\underline a^\epsilon := \underline a \vee \epsilon I$
so that
$$
\bar G^\epsilon(\gamma) := \sup\{\ \frac12 \tr [a \gamma]\ |\
\underline a^\epsilon \le a \le \overline a\ \}.
$$
We then mollify $\bar G^\epsilon$.  Indeed,
let $\eta: \dbS^d \to
[0,1]$ be a regular bump function, i.e., support of $\eta$ is
the unitary ball $O_1$ and $\int_{O_1} \eta(\g) d\g =1$.  We then define
$$
G^\epsilon(\gamma): = \int_{O_1}\
\bar G^\epsilon(\gamma+ \epsilon \gamma^\prime) \
\eta(\gamma^\prime)\  d\gamma^\prime.
$$
It can be shown that
$$
\frac12 \tr[ \underline  a^\epsilon \gamma^\prime]
\le G^\epsilon(\gamma + \gamma^\prime)- G^\epsilon(\gamma)
\le \frac12 \tr[ \overline a \gamma^\prime],
$$
and that there is a constant $C^*$ satisfying
$$
0 \le G^\epsilon(\gamma) - \bar G^\epsilon(\gamma)
\le C^* \epsilon,
$$
where the left inequality thanks to the obvious fact that $\bar G^\e$ is convex.
Moreover $G^\epsilon$ is smooth and convex.
Thus, we can define the Legendre transform of $G^\epsilon$
by
$$
L^\epsilon(a) :=\sup_{\gamma\in \dbS_d}
\{\ \frac12\tr[ a \gamma]- G^\epsilon(\gamma)\ \}.
$$
Then $L^\epsilon(a)$ is finite only if $\underline a^\epsilon \le a \le \overline a$.
Also,  $-C^*  \epsilon \le L^\e(a) \le 0$ for all
$\underline a^\epsilon \le a \le \overline a$ and
$$
G^\epsilon(\gamma) :=\sup_{\underline a^\epsilon \le a \le \overline a}
\{\ \frac12 \tr[ a \gamma]- L^\epsilon(a)\ \}.
$$
We are now ready to prove the approximation result.
Recall that $\cM \subset \dbL_\cP^2$ is the subset
for which the representation \reff{e.main} holds.

\begin{prop}
\label{p.approximate}
Assume that  $\underline{a}$ and $\overline{a}$
satisfy \reff{e.matrix}.   Then, $\cL_{ip} \subset \cM$.
\end{prop}
\proof
Let $\xi \in \cL_{ip}$.  Then $\xi= \f (B_{t_1},\ldots, B_{t_n})$ for
some bounded Lipschitz function $\f$ and $0 \le t_1\le
\ldots \le t_n=1$.  Let $\{v_i\}_{i=1}^n$ be the solutions of \reff{Gpde2}.
Then, $v_i$'s are  bounded and Lipschitz continuous.  Moreover, by the definition
of the $G$-expectations
$$
E^G_t[\xi] = v_i(t,B_{t_1},\ldots, B_{t_{i-1}},B_{t}),
\qquad t \in [t_{i-1},t_{i}).
$$

We approximate $v_i$ as follows.
Let $\f^\epsilon$ be smooth, bounded approximation of
$\f$ so that $\|\f^\epsilon- \f\|_\infty$ tends to zero
and $\|\nabla \f^\epsilon\|_\infty \le \|\nabla \f\|_\infty$.  Define
$v^\epsilon_i(t,x_1,\ldots,x_i,x)$ recursively as in the definition
$G$-expectations in Section \ref{s.G} with data
$\f^\epsilon(B_{t_1},\ldots,B_{t_n})$ and the nonlinearity $G^\epsilon$.
Indeed, $v^\epsilon_i$ is the solution of
\be
\label{e.ge}
-\frac{\partial}{\partial t} v^\epsilon_i(t,x_1,\ldots,x_{i-1},x)
-G^\epsilon(D^2_x v^\epsilon_i(t,x_1,\ldots,x_{i-1},x)) =0,
\ee
on the interval $[t_{i-1},t_{i})$ with final data
$v^\epsilon_i(t_{i},x_1,\ldots,x_{i-1},x)=v^\epsilon_{i+1}(t_{i},x_1,\ldots,x_{i-1},x,x)$.
In the interval $[t_{n-1},1)$, $v^\epsilon_n(t,x_1,\ldots,x_{n-1},x)$ solves
\reff{e.ge} with data $v^\epsilon_n(1,x_1,\ldots,x_{n-1},x)=
{\f}^{\epsilon}(x_1,\ldots,x_{n-1},x)$.

We claim that  the celebrated regularity result of Krylov \cite{krylov} (Theorem 1,
section 6.3, page 292) applies and that
$v^\epsilon_i(t,x_1,\ldots,x_{i-1},x)$ is a smooth function of
$(t,x) \in (t_i,t_{i+1}) \times \dbR^d$.
Indeed, the nonlinearity $G^\epsilon$ depends only on the Hessian
variable.  Moreover, it is constructed so that all its derivatives with respect to 
$\gamma$ are bounded on all of the space.  Hence this nonlinearity $G^\epsilon$
can be directly shown to belong to the class of functions
considered in the Definition 5.5.1 of \cite{krylov}.  Moreover,
in the notation of Theorem 1 of Section 6.3 in \cite{krylov} 
(page 292), the domain $Q= (0,1)\times \dbR^d$.  Therefore,
this theorem
applies to yield existence and interior regularity.  To obtain regularity up
to the terminal condition, we use Theorem 2(b) in \cite{krylov} (Section 6.3,
page 295).  
We may then use the stochastic control representation of this smooth
and classical solution to obtain bounds.  Indeed, the boundedness
and the Lipschitz estimate 
are immediate consequences of the fact that the equation is translation invariant
(or equivalently, the nonlinearity $G^\epsilon$ depends
only on the Hessian).  Hence the solution  is bounded and Lipschitz in all variables.
Moreover the uniform Lipschitz constant of $\f$ is preserved and
for each $i$, we have
\bea
\label{veconv}
\lim_{\e \to 0}\ \|v_i^\epsilon-v_i\|_\infty =0,
\qquad
\sup_{0<\epsilon \le 1}\ \|\nabla v_i^\epsilon\|_\infty \le  \|\nabla \f\|_\infty.
\eea

For $t \in (t_i,t_{i+1})$, we set
\beaa
M^\epsilon_t&:=& v^\epsilon_i(t,B_{t_1},\ldots, B_{t_{i-1}},B_{t}),\\
H^\epsilon_t&:=& \nabla_x v^\epsilon_i(t,B_{t_1},\ldots, B_{t_{i-1}},B_{t}),\\
K^\epsilon_t&:=& G^\epsilon(D^2_x v^\epsilon_i(t,B_{t_1},\ldots, B_{t_{i-1}},B_{t}))
-\frac12 \tr [\hat a_t  D^2_x v^\epsilon_i(t,B_{t_1},\ldots, B_{t_{i-1}},B_{t})],
\eeaa
so that
$$
dM^\epsilon_t = H^\epsilon \cdot dB_t - dK^\epsilon_t.
$$
Let $\cP_\epsilon$ be defined
exactly as $\cP$ but with lower bound $\underline a_\epsilon$
in \reff{abound}.
Then, by the definition of $G^\epsilon$ and $\cP_\epsilon$, we have
that $K^\epsilon$ is non-decreasing $\dbP$ almost surely
for every $\dbP \in \cP_\epsilon$.  But also
since $L^\epsilon \ge - C^* \epsilon$,  we have
\be
\label{e.mart}
-C^*\epsilon \le \sup_{\dbP \in \cP_\epsilon}\ \dbE^\dbP\left[\ -\
K^\epsilon_1\right] \le 0.
\ee

It follows from \reff{veconv} that $M^\epsilon_t$ converges to $M_t:= E^G_t[\xi]$.  Also,
$|H^\epsilon_t|$ is uniformly bounded in $\epsilon$ due to the Lipschitz
estimate on  $v^\epsilon_i$.  Hence $ H^\epsilon \in \cH^2_G$.
Also the Proposition \ref{p.est1} (applied with $\cP_\epsilon$ instead of
$\cP$) yields,
$$
\|K^\epsilon\|_{\dbS^2_{\cP_{\epsilon}}} \le
C \|M^\epsilon\|_{\dbS^2_{\cP_{\epsilon}}} \le C \|\xi\|_\infty.
$$
Moreover, noting that $\cP_\e$ is decreasing as $\e$ increases, 
by Proposition \ref{p.est2}  we obtain the following estimate
$$
\|H^\epsilon - H^{\epsilon^\prime}\|_{\dbH^2_{\cP_{\epsilon_0}}}
+ \|K^\epsilon - K^{\epsilon^\prime}\|_{\dbS^2_{\cP_{\epsilon_0}}}
\le C(\epsilon_0), \qquad 0< \epsilon, \epsilon^\prime \le \epsilon_0,
$$
where
\beaa
C(\epsilon_0)&:=& \sup_{ 0< \epsilon, \epsilon^\prime \le \epsilon_0}\
\left(\|M^\epsilon - M^{\epsilon^\prime}\|_{\dbS^2_{\cP_{\epsilon_0}}}
+\|M^\epsilon - M^{\epsilon^\prime}\|_{\dbS^2_{\cP_{\epsilon_0}}}^{1/2}
\left( \|K^\epsilon\|_{\dbS^2_{\cP_{\epsilon_0}}}+
\|K^{\epsilon^\prime}\|_{\dbS^2_{\cP_{\epsilon_0}}} \right) \right)\\
&\le &  \sup_{ 0< \epsilon, \epsilon^\prime \le \epsilon_0}\
\left(\|M^\epsilon - M^{\epsilon^\prime}\|_{\dbS^2_{\cP_{\epsilon_0}}}
+\|M^\epsilon - M^{\epsilon^\prime}\|_{\dbS^2_{\cP_{\epsilon_0}}}^{1/2}
\left( 2\|\xi\|_\infty \right)\right).
\eeaa
Since $M^\epsilon$ converges uniformly to $M_t$,
$C(\epsilon_0)$ tends to zero with $\epsilon_0$.
Therefore $\{(H^\epsilon, K^\epsilon)\}_{\epsilon}$
is a Cauchy sequence in $\dbH^2_{\cP_{\epsilon_0}}
\times \dbS^2_{\cP_{\epsilon_0}}$ for every
$\epsilon_0$.

By the closure results, Theorem \ref{t.integral} and Theorem \ref{t.Kclose},
we conclude that
there are $H \in  \dbH^2_{\cP_{\epsilon}}$ and $K \in \dbI^2_{\cP_{\epsilon}}$
for every $\epsilon >0$ and that $(M,H,K)$ satisfies
\ref{e.rep} and
\bea
\label{HKuniform}
\|H\|_{\dbS^2_{\cP_\e}} + \|K\|_{\dbS^2_{\cP_\e}} \le C\|\xi\|_{\infty}.
\eea
Clearly $H$ and $K$ are independent of $\e$.  Since by definition and by \reff{abound}
$$
\cP= \cup_{\epsilon >0}\ \cP_\epsilon,
$$
we conclude from the uniform estimates \reff{HKuniform} that
$H \in \dbH^2_{\cP}$, $K \in \dbI^2_{\cP}$.
Moreover, this yields that
$H \in \cH^2_\cP$ and also $-K$
is a $G$-martingale by \reff{e.mart}.
Since $M_t=E^G_t[\xi]$, we have shown that
there is a martingale representation
for the arbitrary random variable $\xi \in \cL_{ip}$.
Hence $\xi \in \cM$.
\ep

\subsection{$\dbL^p_\cP$-spaces}
\label{ss.lp}
In this section we study the properties of the $\dbL^2_\cP$
space. The following result  together
with the example that follows it, imply Lemma  \ref{l.equal}.

\begin{lem}
\label{l.lp1}
For every $p>2$, there exists $C_p$ so that for $\xi \in \cL_{ip}$,
\beaa
\|\xi\|_{\dbL^2_{\cP}} \le C_p\|\xi\|_{\cL^p_G}.
\eeaa
\end{lem}

\proof Since $\xi \in \cL_{ip}$, by its definition in Section 2.1, $M_t := E^G_t[\xi]$ is continuous. 
Moreover, for each $\dbP\in\cP$, by Proposition \ref{p.gcond} we have 
$M_t = \esup_{\dbP' \in \cP(t, \dbP)} \dbE^{\dbP'}_t[\xi]$, $\dbP-$a.s..
Set $M^{*}_t := \sup_{0\le s\le t}M_t$. It suffices to show that
\beaa
\dbE^\dbP[|M^{*}_1|^2] \le C_p\|\xi\|^2_{\dbL^p_G} &\mbox{for all}& \dbP\in\cP.
\eeaa
Now fix $\dbP\in\cP$. Without loss of generality we may assume $\xi\ge 0$.

For any $\l>0$, set $\hat\t := \hat\t_\l := \inf\{t: M_t \ge \l\}$.
Since $M$ is continuous, $\hat\t$ is an $\dbF-$stopping time and
\beaa
\dbP(M^{*}_1 \ge \l) = \dbP(\hat\t \le 1) \le {1\over \l}\dbE^\dbP\Big[M_{\hat\t}\1_{\{\hat\t \le 1\}}\Big].
\eeaa
By Neveu \cite{Neveu} (Proposition VI-1-1), there exist a sequence $\{\dbP_j, j\ge 1\}\subset \cP(\hat\t, \dbP)$
defined in \reff{e.ptp}such that
\beaa
M_{\hat\t} = \sup_{j\ge 1} \dbE^{\dbP_j}_{\hat\t}[\xi],~\dbP-\mbox{a.s.}.
\eeaa
For each $n\ge 1$, denote
\beaa
M^{n}_{\hat\t} := \sup_{1\le j\le n} \dbE^{\dbP_j}_{\hat\t}[\xi].
\eeaa
Then $M^{n}_{\hat\t}\uparrow M_{\hat\t}$, $\dbP-$a.s.. Fix $n$.
Set  $A_j :=
\{M^{n}_{\hat\t}= \dbE^{\dbP_j}_{\hat\t}[\xi]\}$, $1\le j\le n$, and
$\tilde A_1 := A_1$, $\tilde A_j:= A_j \backslash \cup_{1\le i<j} A_i$, $j=2,\cds,n$.
Then $\{\tilde A_j, 1\le j\le n\}\subset \cF^B_{\hat\t}$ form a partition of $\O$.
Define $\hat \dbP^n$ by
\beaa
\hat\dbP^n(E) := \sum_{j=1}^n \dbP_j(E\cap \tilde A_j).
\eeaa
We claim that
\bea
\label{hatPn}
\hat\dbP^n\in \cP({\hat\t},\dbP) &\mbox{and}& M^{n}_{{\hat\t}} 
= \dbE^{\hat\dbP^n}_{\hat\t}[\xi],~~\hat\dbP^n\mbox{a.s.}.
\eea
In fact, $\hat\dbP^n$ is obviously a probability measure and, 
since $\dbP_j \in \cP(\hat\t, \dbP)$,  $\hat\dbP^n = \dbP$ 
on $\cF^B_{\hat\t}$. Then $B$ is a $\hat\dbP$-martingale 
on $[0, \hat\t]$. Moreover, for any stopping time $\t\ge \hat\t$ and any bounded 
$\cF^B_{\t}$-measurable random variable $\eta$, since $B$ is a 
$\dbP_j$-martingale and $\tilde A_j \in \cF^B_{\hat\t}\subset \cF^B_\t$, we have
\beaa
\dbE^{\hat\dbP^n}[B_1 \eta] &=& \sum_{j=1}^n \dbE^{\hat\dbP^n}[B_1 \eta \1_{\tilde A_j}] 
= \sum_{j=1}^n \dbE^{\dbP_j}[B_1 \eta \1_{\tilde A_j}]\\
 &=&  \sum_{j=1}^n \dbE^{\dbP_j}[B_1 \eta \1_{\tilde A_j}] 
 =  \sum_{j=1}^n \dbE^{\dbP_j}[B_\t \eta \1_{\tilde A_j}] = \dbE^{\hat\dbP^n}[B_\t \eta].
\eeaa
Therefore, $\dbE^{\hat\dbP^n}[B_1 |\cF_\t^B] = B_\t$, $\hat\dbP^n$-a.s..
Hence $B$ is a $\hat \dbP^n$-martingale on $[\hat \t, 1]$. So $\hat\dbP^n$ is 
a martingale measure. By \reff{abound}, for each $j$ there exists a constant 
$c_j>0$ so that $BB^T - \overline{a}$ and $BB^T - (c_j I_d\vee \underline{a})$
 are $\dbP_j$-supermartingale and $\dbP_j$-submartingale, respectively. 
Set $\dis c := \min_{1\le j\le n} c_j >0$. Similarly one can show that
 $BB^T - \overline{a}$ and $BB^T - (c I_d\vee \underline{a})$ are 
 $\hat\dbP^n$-supermartingale and $\hat\dbP^n$-submartingale, respectively. 
 This implies that $\hat\dbP^n$ satisfies \reff{abound} and therefore
 $\hat\dbP^n \in \cP({\hat\t},\dbP)$. Finally, for any  bounded $\cF^B_{\hat\t}$-measurable 
 random variable $\eta$, since $\tilde A_j\subset A_j$, we have
\beaa
\dbE^{\hat\dbP^n}[\xi \eta] &=& \sum_{j=1}^n \dbE^{\hat\dbP^n}[\xi \eta \1_{\tilde A_j}] 
= \sum_{j=1}^n \dbE^{\dbP_j}\Big[\dbE^{\dbP_j}_{\hat\t}[\xi] \eta \1_{\tilde A_j}\Big]\\
 &=&  \sum_{j=1}^n \dbE^{\dbP_j}[M_{\hat\t} \eta \1_{\tilde A_j}] 
 = \sum_{j=1}^n \dbE^{\dbP}[M_{\hat\t} \eta \1_{\tilde A_j}] 
 = \dbE^{\dbP}[M_{\hat\t} \eta ]= \dbE^{\hat\dbP^n}[M_{\hat\t} \eta].
\eeaa
Hence $M^{n}_{{\hat\t}} = \dbE^{\hat\dbP^n}_{\hat\t}[\xi]$, $\hat\dbP^n$-a.s.~and
this proves the claim \reff{hatPn}.

Now let $q:=p/(p-1)$ be the conjugate of $p$.
We directly estimate that
\beaa
\dbE^\dbP\Big[M^{n}_{\hat\t}\1_{\{\hat\t \le 1\}}\Big]& =&
\dbE^{\hat\dbP^n}\Big[M^{n}_{\hat\t}\1_{\{\hat\t \le 1\}}\Big]
= \dbE^{\hat\dbP^n}\Big[\dbE^{\hat\dbP^n}_{\hat\t}[\xi]\1_{\{\hat\t \le 1\}}\Big]
= \dbE^{\hat\dbP^n}\Big[\xi\1_{\{\hat\t \le 1\}}\Big]\\
&\le&\Big[\dbE^{\hat\dbP^n}(|\xi|^p)\Big]^{1\over p}\Big[\hat\dbP^n(\hat\t \le 1)\Big]^{1\over q}
=\Big[\dbE^{\hat\dbP^n}(|\xi|^p)\Big]^{1\over p}\Big[\dbP(M^{*}_1 \ge \l)\Big]^{1\over q}\\
&\le& \|\xi\|_{\dbL^p_G}\Big[\dbP(M^{*}_1 \ge \l)\Big]^{1\over q}.
\eeaa
We let $n\to \infty$ to arrive at
\beaa
\dbP(M^{*}_1 \ge \l) \le {1\over \l}\dbE^\dbP\Big[M^{*}_{\hat\t}\1_{\{\hat\t \le 1\}}\Big]
\le \lim_{n \to \infty}\ {1\over \l}\dbE^\dbP\Big[M^{n}_{\hat\t}\1_{\{\hat\t \le 1\}}\Big]
\le {1\over \l} \|\xi\|_{\dbL^p_G}\Big[\dbP(M^{*}_1 \ge \l)\Big]^{1\over q}.
\eeaa
Therefore,
\beaa
\dbP(M^{*}_1 \ge \l) \le {1\over \l^p}\|\xi\|_{\dbL^p_G}^p,
\eeaa
so that for any fixed $\l_0$,
\beaa
\dbE^{\dbP}[|M^{*}_1|^2] &=& 2\int_0^\infty \l \dbP(M^{*}_1\ge \l) d\l \le 2\int_0^{\l_0}\l d\l
+ 2\int_{\l_0}^\infty \l \dbP(M^{*}_T\ge \l) d\l\\
&\le& \l_0^2 + 2\|\xi\|_{\dbL^p_G}^p\int_{\l_0}^\infty{d\l \over \l^{p-1}}
= \l_0^2 + {2\over p-2} \|\xi\|_{\dbL^p_G}^p\l_0^{2-p}.
\eeaa
We choose $\l_0 := \|\xi\|_{\dbL^p_G}$ to conclude that
\beaa
\dbE^{\dbP}[|M^{*}_1|^2] \le C_p\|\xi\|_{\dbL^p_G}^2.
\eeaa
\ep

We next construct a bounded random variable which is not in $\cL^1_G$.
\begin{eg}
\label{e.l1}
{\rm Let $d=1$, $\underline{a}=1$, $\overline{a}=2$,
$E:= \{\limsup_{t\downarrow 0} B_t \slash \sqrt{2t \ln\ln{1\over t}} = 1\}$.
We claim that $\1_E \notin \cL^1_G$.
Indeed, assume that $\1_E \in \cL^1_G$.
Then there exists $\xi_n=\f(B_{t_1},\cds, B_{t_n})\in \cL_{ip}$
such that $\dbE^G[|\xi_n-\1_E|] < {1\over 3}$.
For $\th \in [0,1]$, denote $a^\th_t := 1+\th \1_{[0, t_1)}(t)$
and $\dbP^\th := \dbP^{a^\th}$.
Define $\psi(x) := \dbE^{\dbP_0}[\f(x, x+B_{t_2-t_1},\cds, x+B_{t_n-t_1})]$.
Since $E\in \cF_{0+}\subset \cF_{t_1}$,
for any $\th\in [0,1]$, we have the following
inequality.
\beaa
\dbE^{\dbP^\th}\Big[\Big|\psi(B_{t_1}) - \1_E\Big|\Big]
&=&
\dbE^{\dbP^\th}\Big[\Big|\dbE^{\dbP^\th}_{t_1}
[\f(B_{t_1},\cds, B_{t_n})] - \1_E\Big|\Big]\\
&\le&
\dbE^{\dbP^\th}\Big[\Big|\f(B_{t_1},\cds, B_{t_n}) - \1_E\Big|\Big]
<
{1\over 3}.
\eeaa
Note that $\dbP^0(E) = 1$ and $\dbP^\th(E)=0$ for all $\th>0$. Then
\beaa
\dbE^{\dbP_0}\Big[\Big|\psi(B_{t_1}) - 1\Big|\Big]<{1\over 3} &\mbox{and}
& \dbE^{\dbP^\th}\Big[\Big|\psi(B_{t_1})\Big|\Big] <{1\over 3} ~~\mbox{for all}~\th >0.
\eeaa
The latter implies that
\beaa
\dbE^{\dbP_0}\Big[\Big|\psi(B_{t_1})\Big|\Big] =
\lim_{\th\downarrow 0}\dbE^{\dbP_0}\Big[\Big|
\psi((1+\th)^{1\over 2}B_{t_1})\Big|\Big]
= \lim_{\th\downarrow 0}\dbE^{\dbP^\th}
\Big[\Big|\psi(B_{t_1})\Big|\Big]\le {1\over 3}.
\eeaa
Thus
\beaa
1 \le \dbE^{\dbP_0}\Big[\Big|\psi(B_{t_1}) - 1\Big|\Big]
+ \dbE^{\dbP_0}\Big[\Big|\psi(B_{t_1})\Big|\Big]
\le {1\over 3} + {1\over 3} = {2\over 3},
\eeaa
yielding a contradiction. Hence $1_E \not \in \cL^1_G$.
\ep}\end{eg}

\end{document}